# Asymptotics of the solution of the Cauchy problem for a singularly perturbed system of hyperbolic equations. Part 2. Initial conditions.*

**Nesterov A.V.** 1[0000-0002-4702-4777]

[1] PLEKHANOV Russian University of Economics, Stremyanny lane 36, Moscow, 117997, Russia
`andrenesterov@yandex.ru`

**Abstract.** An asymptotic small parameter expansion of a single Cauchy problem is constructed for a singularly perturbed system of hyperbolic equations describing vibrations of two rigidly connected strings. Equations (such as generalized Korteweg-de Vriesequations) and initial conditions for the terms of the asymptotic expansion of the solution are determined, and under certain assumptions, the residual term is estimated from the residual.

## 1  Introduction

In [1], the author constructed an asymptotic expansion of the solution of a singularly perturbed system of equations of vibrations of interconnected strings

$$\begin{cases} \varepsilon^3(u_{tt} - k_1^2 u_{xx}) = -au + bv + \varepsilon^2 f(u,v), \\ \varepsilon^3(v_{tt} - k_2^2 v_{xx}) = au - bv - \varepsilon^2 f(u,v), \end{cases} \quad (1)$$

with initial conditions

$$\begin{cases} u(x,0) = \overset{0}{u}(\frac{x}{\varepsilon}), u_t(x,0) = \varphi(\frac{x}{\varepsilon}), \\ v(x,0) = \frac{a}{b}\overset{0}{u}(\frac{x}{\varepsilon}), v_t(x,0) = \frac{a}{b}\varphi(\frac{x}{\varepsilon}), \end{cases} \quad (2)$$

where:

$\{u(x,t,\varepsilon), v(x,t,\varepsilon)\}$ - solution, $\{x,t\} \in H = \{|x| < \infty; 0 \le t \le T, T > 0\}$; $0 < \varepsilon \ll 1$.

*The work was carried out with the financial support of the Ministry of Education and Science within the framework of the state assignment in the field of scientific activity on the topic "Models, methods and algorithms of artificial intelligence in economic problems for the analysis and stylization of multidimensional data, time series forecasting and design of recommendation systems", project number FSSW-2023-0004.

small parameter, $k_1, k_2, a > 0, b > 0 - const$, $f(u,v) \subset C^\infty$, in the domain $\Omega = \{|u| < U, |v| < V\}, U, V > 0$, has in $\Omega$ partial derivatives bounded in $\Omega$, $f(0,0) = 0$, the functions $\overset{0}{u}(x), \varphi(x)$ have a certain smoothness, and satisfy the conditions $\overset{0}{u}(x), \overset{0}{v}(x) \in \Omega \forall |x| < \infty$ and fast decreasing at infinity.

*Note 1.* In [1], the system of equations (1) contains a typo that was corrected in this paper: in system (1) $k_1, k_2$, the following words are replaced by $k_1^2, k_2^2$.

*Remark 2.* The initial conditions (2) are consistent in a certain way with respect to functions $\{u(x,t,\varepsilon), v(x,t,\varepsilon)\}$. The approval condition is essential.

Note that the system (1) belongs to the so-called critical case [2] (the degenerate system has a family of solutions). This paper, together with [1], is a continuation of the research initiated in [3], [4].

If the condition of smoothness of the initial conditions is satisfied up to an arbitrary order $N$, the asymptotic expansion (AP) with respect to a small parameter of the solution of problem (1)-(2) in [1] was constructed in the form

$$\begin{pmatrix} u(x,t,\varepsilon) \\ v(x,t,\varepsilon) \end{pmatrix} = \sum_{i=0}^{N} \varepsilon^i \begin{pmatrix} S^I_i u(\zeta_1,t) + S^{II}_i u(\zeta_2,t) \\ S^I_i v(\zeta_1,t) + S^{II}_i v(\zeta_2,t) \end{pmatrix} + \begin{pmatrix} Ru \\ Rv \end{pmatrix} \qquad (3)$$

where $S^I_i u(\zeta_1,t), S^{II}_i u(\zeta_2,t), S^I_i v(\zeta_1,t), S^{II}_i v(\zeta_2,t)$ are terms describing the solution along certain lines - "pseudo-characteristics" of system (1), $R$ is the residual term, stretched variables $\zeta_{1,2}$ have the form $\zeta_{1,2} = (x \mp kt)/\varepsilon$, $\dfrac{bk_1 + ak_2}{a+b} = k, \min(k_1, k_2) < k < \max(k_1, k_2),$

In [1], equations for determining the summands were obtained $S^I_i u(\zeta_1,t), S^{II}_i u(\zeta_2,t), S^I_i v(\zeta_1,t), S^{II}_i v(\zeta_2,t)$, but the question of the initial conditions that these functions must satisfy remained open. This question is answered in this paper.

In view of the fact that in [1] we give a detailed description of obtaining the equations to which the functions obey $S^I_i u(\zeta_1,t), S^{II}_i u(\zeta_2,t), S^I_i v(\zeta_1,t), S^{II}_i v(\zeta_2,t)$, in this paper these calculations are omitted and only the equations themselves are given

$$-S^I_0 u_t + K S^I_0 u_{\zeta_1 \zeta_1 \zeta_1} - h(S_0 u)_{\zeta_1} = 0, \quad S^I_0 v = \frac{a}{b} S^I_0 u,$$

$$S^{II}_0 u_t + K S^{II}_0 u_{\zeta_2 \zeta_2 \zeta_2} - h(S^{II}_0 u)_{\zeta_2} = 0, \quad S^{II}_0 v = \frac{a}{b} S^{II}_0 u, \qquad (4)$$

Here

$$K = \frac{(k^2 - k_2^2)(k^2 - k_1^2)}{2k(a+b)}, h(S_0 u) = -\frac{b}{2k(a+b)}(k^2 - k_2^2) f(S_0 u, \frac{a}{b} S_0 u).$$

The remaining terms *(i>0)* in AP (3) satisfy the linear equations

$$-S^I_i u_t + K S^I_i u_{\zeta_1 \zeta_1 \zeta_1} - (h_S(S^I_0 u) S^I_1 u)_{\zeta_1} = s^I_i u, S^I_i v = \frac{a}{b} S^I_i u + s^I_i v,$$

$$S^{II}_0 u_t + K S^{II}_0 u_{\zeta_2 \zeta_2 \zeta_2} - (h_S(S^{II}_0 u) S^{II}_1 u)_{\zeta_2} = s^{II}_i u, S^{II}_0 v = \frac{a}{b} S^{II}_0 u + s^{II}_i v,$$

where $s^I_i u, s^I_i v, s^{II}_i u, s^{II}_i$ are defined using previously found AP members with numbers *j<i*.

**3 Definition of initial conditions for S-functions**

In [1], the initial conditions for the equations were not defined; the constructed "burst" functions must jointly satisfy the initial conditions (2)

$$S^I u(\zeta_1, t) + S^{II} u(\zeta_2, t)\big|_{t=0} = \overset{0}{u}(x/\varepsilon),$$

$$(S^I u(\zeta_1, t) + S^{II} u(\zeta_2, t))_t\big|_{t=0} = \varphi(x/\varepsilon),$$

$$S^I v(\zeta_1, t) + S^{II} v(\zeta_2, t)\big|_{t=0} = \frac{a}{b} \overset{0}{u}(x/\varepsilon),$$

$$(S^I v(\zeta_1, t) + S^{II} v(\zeta_2, t))_t\big|_{t=0} = \frac{a}{b} \varphi(x/\varepsilon),$$

To determine the initial conditions, we substitute decompositions of functions (3)

$$S^I(\zeta_1, t, \varepsilon) + S^{II}(\zeta_2, t, \varepsilon) = \sum_{i=0}^{\infty} \varepsilon^i \begin{pmatrix} S^I_i u(\zeta_1, t) + S^{II}_i u(\zeta_2, t) \\ S^I_i v(\zeta_1, t) + S^{II}_i v(\zeta_2, t) \end{pmatrix}$$

and derivatives

$$(S^I(\zeta_1, t, \varepsilon) + S^{II}(\zeta_2, t, \varepsilon))_t = \sum_{i=0}^{\infty} \varepsilon^i \begin{pmatrix} S^I_i u(\zeta_1, t) + S^{II}_i u(\zeta_2, t) \\ S^I_i v(\zeta_1, t) + S^{II}_i v(\zeta_2, t) \end{pmatrix}_t$$

in conditions (2). We convert the time derivatives, taking into account the type of arguments $\zeta_{1,2} = (x \mp kt)/\varepsilon$:

$$(S^I(\zeta_1,t,\varepsilon)+S^{II}(\zeta_2,t,\varepsilon))_t = \sum_{i=0}^{\infty}\varepsilon^i \begin{pmatrix} S^I_i u(\zeta_1,t)+S^{II}_i u(\zeta_2,t) \\ S^I_i v(\zeta_1,t)+S^{II}_i v(\zeta_2,t) \end{pmatrix}_t =$$

$$= \sum_{i=0}^{\infty}\varepsilon^i \begin{pmatrix} \dfrac{-k}{\varepsilon}(S^I_i u)_{\zeta_1}+(S^I_i u)_t+\dfrac{k}{\varepsilon}(S^{II}_i u)_{\zeta_2}+(S^{II}_i u)_t \\ \dfrac{-k}{\varepsilon}(S^I_i v)_{\zeta_1}+(S^I_i v)_t+\dfrac{k}{\varepsilon}(S^{II}_i v)_{\zeta_2}+(S^{II}_i v)_t \end{pmatrix}$$

By standard procedure, we obtain the initial conditions for the functions $S^I_i u(x/\varepsilon,0), S^{II}_i u(x/\varepsilon,0), S^I_i v(x/\varepsilon,0), S^{II}_i v(x/\varepsilon,0), i=0,1,...$

$$\sum_{i=0}^{\infty}\varepsilon^i \begin{pmatrix} S^I_i u(\zeta_1,0)+S^{II}_i u(\zeta_2,0) \\ S^I_i v(\zeta_1,0)+S^{II}_i v(\zeta_2,0) \end{pmatrix} = \begin{pmatrix} \overset{0}{u}(x/\varepsilon) \\ \dfrac{a}{b}\overset{0}{u}(x/\varepsilon) \end{pmatrix},$$

$$\sum_{i=0}^{\infty}\varepsilon^i \begin{pmatrix} \dfrac{-k}{\varepsilon}(S^I_i u)_{\zeta_1}+(S^I_i u)_t+\dfrac{k}{\varepsilon}(S^{II}_i u)_{\zeta_2}+(S^{II}_i u)_t \\ \dfrac{-k}{\varepsilon}(S^I_i v)_{\zeta_1}+(S^I_i v)_t+\dfrac{k}{\varepsilon}(S^{II}_i v)_{\zeta_2}+(S^{II}_i v)_t \end{pmatrix} = \begin{pmatrix} \varphi(x/\varepsilon) \\ \dfrac{a}{b}\varphi(x/\varepsilon) \end{pmatrix}.$$

(5)

We restrict ourselves to obtaining initial conditions only for the principal members of the AP $S^I_0 u(x/\varepsilon,0), S^{II}_0 u(x/\varepsilon,0)$. Equating the terms with the principal powers of the parameter $\varepsilon$ in (5), we obtain

$$\begin{pmatrix} S^I_0 u(\zeta_1,0)+S^{II}_0 u(\zeta_2,0) \\ S^I_0 v(\zeta_1,0)+S^{II}_0 v(\zeta_2,0) \end{pmatrix} = \begin{pmatrix} \overset{0}{u}(x/\varepsilon) \\ \dfrac{a}{b}\overset{0}{u}(x/\varepsilon) \end{pmatrix} \qquad (6)$$

$$\begin{pmatrix} \dfrac{-k}{\varepsilon}(S^I_0 u)_{\zeta_1}+(S^I_i u)_t+\dfrac{k}{\varepsilon}(S^{II}_0 u)_{\zeta_2}+(S^{II}_0 u)_t \\ \dfrac{-k}{\varepsilon}(S^I_0 v)_{\zeta_1}+(S^I_0 v)_t+\dfrac{k}{\varepsilon}(S^{II}_i v)_{\zeta_2}+(S^{II}_0 v)_t \end{pmatrix} = \begin{pmatrix} \varphi(x/\varepsilon) \\ \dfrac{a}{b}\varphi(x/\varepsilon) \end{pmatrix}$$

Since $S^J_0 v = \dfrac{a}{b} S^J_0 u, J=I,II$ then, when the conditions for $S^J_0 u, J=I,II$ are met, the conditions for $S^J_0 v, J=I,II$ are met, so only conditions on can be considered $S^J_0 u, J=I,II$. Leaving the main terms in (6), we obtain

$$S^I_{\ 0}u(\zeta_1,0) + S^{II}_{\ 0}u(\zeta_2,0) = \overset{0}{u}(x/\varepsilon),$$
$$\left. \frac{-k}{\varepsilon}(S^I_{\ 0}u)_{\zeta_1} + \frac{k}{\varepsilon}(S^{II}_{\ i}u)_{\zeta_2} \right|_{t=0} = 0. \tag{7}$$

Let's get rid of the small parameter in the second condition (7)

$$S^I_{\ 0}u(x/\varepsilon,0) + S^{II}_{\ 0}u(x/\varepsilon,0) = \overset{0}{u}(x/\varepsilon),$$
$$-(S^I_{\ 0}u)_{\zeta_1} + (S^{II}_{\ 0}u)_{\zeta_2} = 0 \tag{8}$$

If conditions (8) are met, the initial conditions for the function $S^I_{\ 0}v(x/\varepsilon,0), S^{II}_{\ 0}v(x/\varepsilon,0)$ are met.

We denote $x/\varepsilon = \xi$ and rewrite the initial conditions (8). Note that
$$Su(\zeta,t)_t = Su(\zeta(x,t),t)_t = Su(\zeta,t)_t + Su(\zeta,t)_\zeta \zeta_t =$$
$$= Su(\zeta,t)_t + \varepsilon^{-1}(\pm)kSu(\zeta,t)_\zeta$$

From the first condition (8), given that
$$\zeta_1|_{t=0} = \zeta_2|_{t=0} = \xi = x/\varepsilon,$$

we get
$$S^I_{\ 0}u(\xi,0) + S^{II}_{\ 0}u(\xi,0) = \overset{0}{u}(\xi), \tag{9}$$

We convert the second condition (8) to derivatives:

$$(-kS^I_{0\ \zeta}u(\zeta_1,t) + kS^{II}_{0\ \zeta}u(\zeta_2,t))\Big|_{t=0} = k(-S^I_{0\ \zeta}u(\zeta_1,t) + S^{II}_{0\ \zeta}u(\zeta_2,t))\Big|_{t=0} =$$
$$= k(-S^I_{0\ \xi}u(\xi,0) + S^{II}_{0\ \xi}u(\xi,0)) = k(-S^I_0 u(\xi,0) + S^{II}_0 u(\xi,0))_\xi = 0,$$

From here
$$(-S^I_0 u(\xi,0) + S^{II}_0 u(\xi,0))_\xi = 0$$

Therefore, the sum on the left side of the equation is a constant $-S^I_0 u(\xi,0) + S^{II}_0 u(\xi,0) = C$.

Given the requirement
$$\{S^I_0 u(\xi,0) \to 0, S^{II}_0 u(\xi,0) \to 0; |\xi| \to \infty\}$$

we get that *C=0* and
$$S^I_0 u(\xi,0) = S^{II}_0 u(\xi,0) \tag{10}$$

From (9) and (10) we obtain
$$S^I_{\ 0}u(\xi,0) = S^{II}_{\ 0}u(\xi,0) = \frac{1}{2}\overset{0}{u}(\xi), \tag{11}$$

For the following terms of expansion (3), the initial conditions are obtained similarly.

Accordingly $S^I{}_0 u(\xi,0), S^{II}{}_0 u(\xi,0)$, there are solutions of initial problems for generalized Korteweg-de Vries equations

$$-S^I{}_0 u_t + K S^I{}_0 u_{\zeta_1 \zeta_1 \zeta_1} - h(S_0 u)_{\zeta_1} = 0, |\zeta_1| < \infty, t > 0,$$
$$S^I{}_0 u(\zeta_1, 0) = \frac{1}{2} \overset{0}{u}(\xi), \quad (12)$$

$$S^{II}{}_0 u_t + K S^{II}{}_0 u_{\zeta_2 \zeta_2 \zeta_2} - h(S^{II}{}_0 u)_{\zeta_2} = 0, |\zeta_2| < \infty, t > 0,$$
$$S^{II}{}_0 u(\zeta_2, 0) = \frac{1}{2} \overset{0}{u}(\xi). \quad (13)$$

The remaining terms of expansion (3) are solutions of Cauchy problems for equations (5).

***Estimates of S functions***

Impose condition I.

Condition I.

Let the initial conditions and function $h(z)$ in problems (12), (13) be such that there exists $T>0$, such that $\forall\ 0 \leq t \leq T, |x| < \infty$ the solution of problems (12), (13) exists, is unique, and bounds on solutions and derivatives up to the third order are satisfied.

$$\left|(S^I{}_0 u)_\zeta^{(k)}\right| < C\theta(\zeta_1), \left|(S^{II}{}_0 u)_\zeta^{(k)}\right| < \theta(\zeta_2), \theta(\zeta) \underset{|\zeta| \to \infty}{\to} 0. \quad (14)$$

*Remark.* It is desirable to have verifiable conditions (I) for the data of equations (12), (13) and initial data, as well as to have estimates of the rate of decrease of solutions to problems (12), (13) with rapidly decreasing initial conditions. Numerous papers devoted to the study of the KdF equation (for example, papers [5]-[9]) mainly study soliton and soliton-like solutions, as well as the asymptotics of the solution for t→∞.

***Evaluation of the residual term.***

The final constructed AP has the form

$$\begin{pmatrix} u(x,t,\varepsilon) \\ v(x,t,\varepsilon) \end{pmatrix} = \begin{pmatrix} S^I u(\zeta_1,t,\varepsilon) \\ S^I v(\zeta_1,t,\varepsilon) \end{pmatrix} + \begin{pmatrix} S^{II} u(\zeta_2,t,\varepsilon) \\ S^{II} v(\zeta_2,t,\varepsilon) \end{pmatrix} + R =$$
$$= \sum_{i=0}^N \varepsilon^i \left( \begin{pmatrix} S^I u(\zeta_1,t) \\ S^I v(\zeta_1,t) \end{pmatrix} + \begin{pmatrix} S^{II} u(\zeta_2,t) \\ S^{II} v(\zeta_2,t) \end{pmatrix} \right) + R = \quad (15)$$
$$= S^I(\zeta_1,t,\varepsilon) + S^{II}(\zeta_2,t,\varepsilon) + R = U + R$$

We give an estimate of the residual term with respect to the residual for the case $N=0$

$$\begin{pmatrix} u(x,t,\varepsilon) \\ v(x,t,\varepsilon) \end{pmatrix} = \left( \begin{pmatrix} S^I{}_0 u(\zeta_1,t) \\ S^I{}_0 v(\zeta_1,t) \end{pmatrix} + \begin{pmatrix} S^{II}{}_0 u(\zeta_2,t) \\ S^{II}{}_0 v(\zeta_2,t) \end{pmatrix} \right) + R$$

In accordance with the AR construction algorithm, the residual terms $Ru, Rv$ satisfy the problem

$$\begin{aligned} \varepsilon^3 (Ru_{tt} - k_1^2 Ru_{xx}) &= -(aRu - bRv) + \varepsilon^2 Rf(Ru, Rv) + ru, \\ \varepsilon^3 (Rv_{0,tt} - k_2^2 Rv_{xx}) &= (aRu - bRv) - \varepsilon^2 Rf(Ru, Rv) + rv, \end{aligned} \qquad (16)$$

and initial conditions

$$\begin{cases} Ru(x,0) = 0, Ru_t(x,0) = 0, \\ Rv(x,0) = 0, Rv_t(x,0) = 0, \end{cases} \qquad (17)$$

Here's the $Rf$ function

$$Rf(Ru, Rv) = (f(S^I u + S^{II} u + Ru, S^I v + S^{II} v + Rv) - \\ -(f(S^I u + S^{II} u, S^I v + S^{II} v), \\ Rf(0,0) = 0$$

From the estimates (14), the estimate follows

$$ru = O(\varepsilon(\theta(\zeta_1) + \theta(\zeta_2))), \\ rv = O(\varepsilon(\theta(\zeta_1) + \theta(\zeta_2))).$$

*Results.*

1. Initial conditions are obtained for determining the principal terms of the AP solution of problem (1)-(2) - the Cauchy problem with initial conditions of the "narrow cap" type for a singularly perturbed system of hyperbolic equations (1) in the critical case.
2. If condition (I) is satisfied, an estimate is made for the residual term discrepancy in the AP.

Questions remain open.

1. Conditions for the existence and evaluation of functions $S$..
2. Estimates of the residual term in the energy norm.

**Conclusion.**
1. The system of equations (1) can be considered as one of the variants of generalization of the problem of coupled oscillators, considered in many works, in particular, in [5],[9], in which the phenomenon of the Pasta-Ulam-Fermi paradox is analyzed and the connection of the system of coupled oscillators with the Korteweg - de Vries equation is noted.

2. In this formulation, the AR solutions of system (1) with initial conditions of the burst type (2) in the first approximation have the form (3), where $S^I_0 u, S^{II}_0 u$ they are described by the initial problems (12), (13) for generalized Korteweg-de Vries equations.

3. If the initial conditions for the functions are inconsistent $u, v$, the solution will contain components of the type of undamped oscillations with an asymptotically large frequency.

4. In the future, we plan to study the asymptotics of solving the Cauchy problem for system (1) with the addition of dissipative terms.